\documentclass{article}
\usepackage{url,amsthm,mathrsfs,amsmath}
\theoremstyle{plain}  
\newtheorem{theorem}{Theorem}

\def\D{{\nabla_\infty}}
\def\Dbar{{{\overline\nabla}_\infty}}
\usepackage{amsmath,mathrsfs}
\RequirePackage[colorlinks,citecolor=blue,urlcolor=blue]{hyperref}      
\makeatletter 
\def\@biblabel#1{\hspace*{-\labelsep}}
\makeatother
\allowdisplaybreaks
\begin{document}
\title{A property of Petrov's diffusion}
\author{S. N. Ethier\thanks{Department of Mathematics, University of Utah, 155 South 1400 East, Salt Lake City, UT 84112, USA. e-mail: ethier@math.utah.edu.  Partially supported by a grant from the Simons Foundation (209632).}}
\date{}
\maketitle

\begin{abstract}
Petrov constructed a diffusion process in the Kingman simplex whose unique stationary distribution is the two-parameter Poisson--Dirichlet distribution of Pitman and Yor.  We show that the subset of the simplex comprising vectors whose coordinates do not sum to 1 acts like an entrance boundary for the diffusion.\medskip\par
\end{abstract}

\section{Introduction}

Petrov (2009) constructed an infinite-dimensional diffusion process in the compact Kingman simplex
$$
\Dbar:=\bigg\{x = (x_1,x_2,\ldots): x_1 \ge x_2 \ge \cdots \ge0,\; \sum_{i=1}^\infty x_i \le 1\bigg\}
$$
depending on two parameters, $\alpha$ and $\theta$ with $0\le\alpha<1$ and $\theta>-\alpha$. Its generator is\footnote{Petrov omitted the common factor of $\frac{1}{2}$ for simplicity but we include it so that formulas are consistent with those in the literature.}
$$
A:=\frac{1}{2}\sum_{i,j=1}^\infty x_i(\delta_{ij}-x_j)\frac{\partial^2}{\partial x_i\,\partial x_j}
-\frac{1}{2}\sum_{i=1}^\infty(\theta x_i+\alpha)\frac{\partial}{\partial x_i}
$$
acting on the subalgebra of $C(\Dbar)$ generated by the sequence of functions $\varphi_1,\varphi_2,\varphi_3,\ldots$ defined by 
$$
\varphi_m(x):=\sum_{i=1}^\infty x_i^m,\quad m=2,3,\ldots,\qquad \varphi_1(x):=1.
$$
More precisely, for $\varphi\in\mathscr{D}(A)$, $A\varphi$ is evaluated on the dense subset
$$
\D:=\bigg\{x = (x_1,x_2,\ldots): x_1 \ge x_2 \ge \cdots \ge0,\; \sum_{i=1}^\infty x_i = 1\bigg\}
$$
and extended to $\Dbar$ by continuity.  For example, 
$$
A\varphi_m=\frac{1}{2}m[(m-1-\alpha)\varphi_{m-1}-(m-1+\theta)\varphi_m],\quad m=2,3,\ldots.
$$
The unique stationary distribution of Petrov's diffusion is Pitman and Yor's (1997) two-parameter generalization of the Poisson--Dirichlet distribution, which we denote by $\text{PD}_{\alpha,\theta}$ and regard as a Borel probability measure on $\Dbar$ that is concentrated on $\D$.

The special case $\alpha=0$ (and hence $\theta>0$) is the unlabeled infinitely-many-neutral-alleles diffusion model of population genetics; see Ethier and Kurtz (1981).  Its unique stationary distribution is of course the (one-parameter) Poisson--Dirichlet distribution, $\text{PD}_{0,\theta}$, of Kingman (1975).  Petrov's diffusion was not motivated by population genetics; however, see De Blasi, Ruggiero, and Span\`o (2014).

Feng, Sun, Wang, and Xu (2011) derived an explicit formula for the transition density $p(t,x,y)$ of Petrov's diffusion with respect to $\text{PD}_{\alpha,\theta}$, which had earlier been done in the special case $\alpha=0$ by Ethier (1992).  Recently, Zhou (2013) found an elegant simplification of this formula.  However, all we will need here is the fact that the Feller semigroup $\{T(t)\}$ on $C(\Dbar)$ generated by the closure of $A$ can be expressed as
\begin{equation}\label{density}
T(t)f(x)=\int_\Dbar f(y)\,p(t,x,y)\,\text{PD}_{\alpha,\theta}(dy),\quad x\in\Dbar,\,t>0,
\end{equation}
for all $f\in C(\Dbar)$.

In particular, letting $\{X_t,\,t\ge0\}$ denote Petrov's diffusion, it follows that $P_x(X_t\in\D)=1$ for every $x\in\Dbar$ and $t>0$, where the subscript $x$ denotes the initial state.  A question left open by Petrov (2009) is whether the stronger statement,
\begin{equation}\label{entrance}
P_x(X_t\in\D\text{ for all }t>0)=1,\quad x\in \Dbar,
\end{equation}
holds.  

In the special case $\alpha=0$, (\ref{entrance}) was proved by Ethier and Kurtz (1981).  It was later realized that this result has a simple interpretation.  The unlabeled infinitely-many-neutral-alleles diffusion model has a more informative labeled version, namely the Fleming--Viot process in $\mathscr{P}(S)$ (the set of Borel probability measures on the compact set $S$ with the topology of weak convergence) with mutation operator
$$
Bg(z):=\frac{1}{2}\theta\int_S(g(\zeta)-g(z))\,\nu_0(d\zeta),
$$
where $\nu_0\in\mathscr{P}(S)$ is nonatomic.  The unlabeled model is a transformation of the labeled one.  The transformation takes $\mu\in\mathscr{P}(S)$ to $x\in\Dbar$, where $x$ is the vector of descending order statistics of the sizes of the atoms of $\mu$.  Then (\ref{entrance}) is equivalent to the assertion that the Fleming--Viot process, regardless of its initial state, instantly becomes purely atomic and remains so forever.

If $\alpha>0$, there is no such interpretation of (\ref{entrance}) because whether there is a Fleming--Viot  process corresponding to Petrov's diffusion is unknown.  In fact, this is an open problem that was posed by Feng (2010, p.~112).

Notice that another way to express (\ref{entrance}) is to say that $\Dbar-\D$ acts like an entrance boundary for the diffusion.  Technically, this is not quite accurate because $\D$ has no interior, so its boundary is all of $\Dbar$.  But what we mean is simply that, starting at a state $x\in\Dbar-\D$, the process instantly enters $\D$ and never exits.

A weaker version of (\ref{entrance}) was obtained by Feng and Sun (2010) using the theory of Dirichlet forms.  They showed that
\begin{equation}\label{stationary}
\int_\Dbar P_y(X_t\in\D\text{ for all }t>0)\,\text{PD}_{\alpha,\theta}(dy)=1,
\end{equation}
which is to say that the stationary version of Petrov's diffusion has $\D$ as its natural state space.

In the next section we will use (\ref{density}) and (\ref{stationary}) to prove (\ref{entrance}).

\section{Entrance boundary}

Let us begin by showing why Ethier and Kurtz's (1981) proof of (\ref{entrance}) when $\alpha=0$ fails when $\alpha>0$.  

First, we extend the domain of $A$.  Let $\mathscr{H}:=\{h\in C^2[0,1]: h(0)=h'(0)=0\}$, and for $h\in\mathscr{H}$ define $\psi_h\in C(\Dbar)$ by
$$
\psi_h(x):=\sum_{i=1}^\infty h(x_i).
$$
Let $A^+$ be $A$ acting on the subalgebra of $C(\Dbar)$ generated by $\{1\}\cup\{\psi_h: h\in\mathscr{H}\}$.  Again, for $\varphi\in\mathscr{D}(A^+)$, $A^+\varphi$ is evaluated on $\D$ and extended to $\Dbar$ by continuity.  It is easy to see that $A^+\subset \overline{A}$.  We also notice that $\varphi_m\in\mathscr{D}(A^+)$ for all real $m\ge2$ (not just integers), where
$$
\varphi_m(x):=\sum_{i=1}^\infty x_i^m,\quad m\in(0,1)\cup(1,\infty),\qquad \varphi_1(x):=1.
$$ 
This leads to the conclusion that, for every real $m\ge2$,
\begin{eqnarray*}
Z_m(t)&:=&\varphi_m(X_t)-\varphi_m(X_0)\\
&&\quad{}-\frac{1}{2}m\int_0^t[(m-1-\alpha)\varphi_{m-1}(X_s)-(m-1+\theta)\varphi_m(X_s)]\,ds
\end{eqnarray*}
is a continuous square-integrable martingale with increasing process
$$
I_m(t)=m^2\int_0^t (\varphi_{2m-1}-\varphi_m^2)(X_s)\,ds.
$$
A difficulty occurs when trying to extend this last conclusion to $1<m<2$.  Fix such an $m$ and define $h_\varepsilon\in\mathscr{H}$ for $\varepsilon>0$ by
$h_\varepsilon(u):=(u+\varepsilon)^m -\varepsilon^m-m\varepsilon^{m-1}u$. 
Then $A^+\psi_{h_\varepsilon}$ converges pointwise as $\varepsilon\to0+$ but not boundedly because of two awkward terms:
$$
\frac{1}{2}\sum_{i=1}^\infty x_i h_\varepsilon''(x_i)\quad\text{and}\quad
-\frac{1}{2}\alpha\sum_{i=1}^\infty h_\varepsilon'(x_i).
$$
Both sums are monotonically increasing as $\varepsilon$ decreases to 0, but their coefficients have opposite signs, so the monotone convergence theorem does not apply if $\alpha>0$, and we cannot show that $\int_0^t\varphi_{m-1}(X_s)\,ds\in L^2$ for all $t>0$.

We are now ready for the proof of (\ref{entrance}).

\begin{theorem}
Eq.\ (\ref{entrance}) holds for Petrov's diffusion.  
\end{theorem}

\begin{proof}
It is enough to show that 
$$
P_x(X_t\in\D\text{ for all }t\ge s)=1
$$ 
for every $x\in \Dbar$ and $s>0$.  By (\ref{stationary}),
\begin{equation}\label{ae}
P_y(X_t\in\D\text{ for all }t\ge0)=1\;\; \text{a.e.-PD}_{\alpha,\theta}(dy).
\end{equation}
Let $s>0$ be arbitrary.  Then, for every $x\in\Dbar$,
\begin{eqnarray*}
&&\!\!\!\!\!\!\!\!\!\!\!\!P_x(X_t\in\D\text{ for all }t\ge s)\\
&=&E_x[P_x(X_t\in\D\text{ for all }t\ge s\mid X_r,\,0\le r\le s)]\\
&=&E_x[P_{X_s}(X_t\in\D\text{ for all }t\ge0)]\\
&=&\int_\Dbar P_y(X_t\in\D\text{ for all }t\ge0)\,p(s,x,y)\,\text{PD}_{\alpha,\theta}(dy)\\
&=&\int_\Dbar p(s,x,y)\,\text{PD}_{\alpha,\theta}(dy)\\
&=&P_x(X_s\in\Dbar)\\
&=&1,
\end{eqnarray*}
where the third equality uses (\ref{density}) and the fourth equality uses (\ref{ae}).
\end{proof}

\end{document}